
\magnification=\magstep1
\baselineskip=14pt
\def \Box {\hbox{}\nobreak
 \vrule width 1.6mm height 1.6mm depth 0mm  \par \goodbreak \smallskip}

\vskip .4cm \centerline{\bf FINDING SPARSE SYSTEMS OF PARAMETERS}
\vskip .3cm

\centerline{\bf David Eisenbud}
\centerline{\sevenrm Brandeis University, Waltham MA 02254}
\centerline{\sevenrm eisenbud@math.brandeis.edu} \vskip .1cm
\centerline{and} \vskip .1cm
\centerline{{\bf Bernd Sturmfels}}
\centerline{\sevenrm Cornell University, Ithaca, NY 14853}
\centerline{\sevenrm bernd@math.cornell.edu}
\vskip .4cm

\midinsert \narrower
\noindent {\bf Abstract:  }
{For several computational procedures such as finding
radicals and Noether normalizations, it is important
to choose as sparse as possible a system of parameters in
a polynomial ideal or modulo a polynomial ideal.
We describe new strategies for these tasks, thus providing
solutions to problems (1) and (2) posed in
[Eisenbud-Huneke-Vasconcelos 1992].

To accomplish the first task we introduce a notion of ``setwise
complete intersection".  We prove that a set of monomials
generating an ideal of codimension $c$ in a polynomial ring can be
partitioned into $c$ disjoint sets forming a setwise complete
intersection, although the corresponding result is false for
arbitrary sets of polynomials. We reduce the general case to
the monomial case by a deformation argument.
For homogeneous ideals the output
is homogeneous. Our analysis of the second task is based on a
concept of Noether complexity for homogeneous ideals and its
characterization in terms of Chow forms. }
 \endinsert

\vskip 1cm

\centerline{\bf Introduction}

\vskip .2cm

\noindent  Let $k$ be a field and let $S := k[x_1,\ldots,x_m]$
be the polynomial
ring. Let $J$ be an ideal of $S$, possibly 0, and let $R = S/J$.
Given a finite set ${\cal F} \subset R$, generating a proper
ideal $I$, it is a prerequisite for many algebraic computations
to  find a maximal system of parameters in the ideal $I$.
By this we mean a system of
$\,codim  (I)$ elements of $I$ which generate an ideal of the
same codimension; see for example
[Eisenbud-Huneke-Vasconcelos 1992,  Eisenbud 1993,  Krick-Logar 1991,
and Vasconcelos 1993a,b].
The feasibility of the subsequent computation often
hinges on the fact that the system of parameters is ``nice'';
typically, that it consists of polynomials which are
 reasonably sparse, and of low degree.

In this paper we address the question of  how to compute
a system of parameters which is  as sparse as possible.
Given a set of polynomials $\cal F$ that generate an
ideal of codimension $c$, we study the ways of dividing $\cal F$
into subsets
$\, {\cal F}_1,{\cal F}_2, \ldots,{\cal F}_c \,\subseteq \,{\cal F} \,$
such that if $f_i$ is a sufficiently general linear combination of
elements from ${\cal F}_i$ then $f_1,\dots, f_c$ form a system
of parameters modulo $J$.

In the first  section we treat the case $J=0$.
This arises when computing the radical of an ideal.
Our main result in this case is Theorem 1.3, which says
 that if the elements of $\cal F$ are monomials, and somewhat
more generally, then
the ${\cal F}_i$ may be chosen to be disjoint subsets of $\cal F$. We
give examples to show that this result fails for arbitrary polynomials,
but we can reduce the general case to this one by a deformation
argument, using partial Gr\"obner bases.  Although it is unpleasant
to have to compute a Gr\"obner basis for this purpose, and the
worst-case complexity
of the computation is certainly very bad, the payoff is high:  if the
input polynomials are homogeneous of varying degrees, the output
can be made homogeneous without any loss of sparseness from the
inhomogeneous case.

In the second section we consider the case of arbitrary $J$,
and we present a simple ``greedy'' algorithm (2.1). We then concentrate
on the important case when $J$ is homogeneous and unmixed,
and $I$ is the ideal $(x_1,\ldots,x_m)$. This is the
problem of Noether normalization.
In this case we
compare some possible meanings of the term ``sparse'', with
the conclusion that the ``correct'' measure of sparseness will
vary with the application at hand.  Choosing one of these,
we define the Noether complexity, which is
the sparseness of the sparsest possible Noether normalization.
It is characterized in terms of the Chow form of $J \,$ (Theorem 2.7)
and can thus be computed in single-exponential time in $m$
(cf.~[Caniglia 1990]). Another simple approach to computing a
Noether normalization of $J$ is to lift one from any initial
ideal of $J$; this lifting method is explained following Proposition
2.8. We demonstrate by explicit examples that, in general,
neither this lifting method nor the greedy algorithm
attains the Noether complexity.

To simplify the discussion, we assume throughout that $k$ is
an infinite field, though in practice any sufficiently large field
will do. In order to use our algorithms it is necessary to compute
the codimensions of various ideals.  Good methods for doing
this are discussed in Bayer-Stillman [1992] and
Bigatti-Caboara-Robbiano [1993].

\vskip .1cm

Both authors were partially supported by the  NSF during the
preparation of this paper. The second author was also supported
by an A.P.~Sloan Fellowship.

\vskip .5cm

\beginsection 1.  Systems of parameters in a polynomial ring

We retain the notation of the introduction.
In this section we treat the case $J=0$, and assume that
${\cal F}$ is a subset of a
proper ideal $ I $ of the polynomial ring $S := k[x_1,\ldots,x_m]$.

\vskip .3cm

\centerline{\bf   Some obvious approaches  and their drawbacks}

\vskip .1cm

\noindent Let $c$ be the codimension of $I$.  A set of $c$ linear
combinations of the polynomials in ${\cal F}$ with  sufficiently
general coefficients in $k$ is a system of parameters, but unfortunately
does
not have the desired sparseness.
Further, if the polynomials in ${\cal F}$ are homogeneous but of
different degrees and a homogeneous system of parameters is
required, then in this approach one must first  replace ${\cal
F}$ by a set of polynomials all of the same degree, for example by
multiplying each one by a power of a generic linear form, or by
replacing each by the ideal it generates in degree equal to the
maximal degree in ${\cal F}$.  This process
dramatically destroys sparseness, and raises the degrees of
the elements of ${\cal F}$ in a way that seems unnecessary.

It is thus natural to ask for the  smallest subsets
$\, {\cal F}_1,{\cal F}_2, \ldots,{\cal F}_c \,\subseteq \,{\cal F} \,$
such that the linear combinations
$$ f_1 \,\,= \sum_{f\in {\cal F}_1} \,r_{1,f} \cdot f \, , \quad
 f_2 \,\,= \sum_{f\in {\cal F}_2} \,r_{2,f} \cdot f \, , \,\,
\ldots \, , \quad f_c \,= \,\sum_{f\in {\cal F}_c} \,r_{c,f} \cdot f
\eqno (*) $$
generate an ideal of codimension $c$ (that is, form a system of
parameters) for some choice of coefficients $r_{i,f}$.
Supposing that no proper subset of ${\cal F}$ generates
an ideal of codimension $c$, an optimal result of this type would be
to take  ${\cal F}_1,{\cal F}_2,\dots,{\cal F}_c $ to be a
partition of ${\cal F}$, that is, disjoint subsets whose union is ${\cal F}$.

A first hope might be that one could define the sets ${\cal F}_j$
inductively by the condition that $${\cal F}_1 \, \cup \,{\cal F}_2 \,
\cup \, \ldots \cup \,{\cal F}_j$$ is the smallest initial subset generating
an ideal of codimension $j$. But this is wrong even for monomial
ideals as the following  example from [Eisenbud 1993] shows:

\vskip .2cm

\noindent {\bf Cautionary Example 1.1.}  Let $m=4$ and
$ \,{\cal F} \,:=\, \{x_1 x_2 , x_2 x_3, x_4^2, x_1 x_3\}$.
We have  $c=3$ and the partition suggested above is
$$\{x_1 x_2\},\,\{x_2 x_3, x_4^2\},\,\{x_1 x_3\}.$$ However,
no sequence of the form
$$ x_1 x_2 \,, \,\lambda x_2 x_3 +\mu x_4^2 \,,\,x_1 x_3$$
is a system of parameters, since it is
contained in the height 2 ideal $(x_1,\lambda x_2 x_3 +\mu x_4^2)$.  On the
other hand, the partition $$\{x_1 x_2\},\{x_2 x_3,x_1 x_3\},\{x_4^2\}$$
does have the
desired property. \Box

\vskip .1cm

Unfortunately, partitions with the desired property need not exist.
The following example was worked out in conversation with Joe Harris.

\vskip .3cm

\noindent {\bf Cautionary Example 1.2.  }   Let
$\{C_{ij}\}_{1\leq i \leq j\leq 5}$ be any $10$ distinct
irreducible space curves in ${\bf P}^3$, and take $d$ an integer
large enough so that for each $i$ the ideal of the union of the $6$
$C_{pq}$ whose indices $p,q$ do not include
$i$ is generated by forms of degree $d$.  For $i=1,\dots,5$ let
$g_i$ be a general form of degree $d$ vanishing on these $6$
curves $C_{pq}$.  Let ${\cal F} = \{g_1,\ldots, g_5 \}$.
It is easy to see that $\cal F$ has no zeros in ${\bf P}^3$
and hence generates an ideal of codimension 4.

We claim that there is no partition of ${\cal F}$  into four disjoint subsets
${\cal F}_i$ and choice of coefficients $\,r_{i,f}\,$ such that
$(*)$ generates an ideal of codimension 4. If such
a partition existed, then 3 of the sets ${\cal F}_i$ would have to be
singletons. Hence some 3 forms $g_i,g_j,g_k$ would have to be a
system of parameters, and vanish only at
finitely many points in ${\bf P}^3$.
Since $g_i,g_j,g_k$  vanish on the curve $C_{uv}$,
where  $\{i,j,k,u,v\}=\{1,\dots,5\}$, this is impossible.

There is no example of this type with $d=2$, but here is
one with $d=3$:
 Let  $p_1,\ldots, p_5$ be general points in ${\bf P}^3$,  let
$C_{ij}$ be the line passing through $p_i$ and $p_j$, and let $l_{ijk}$
be the equation  of the plane containing the three points
$p_i,p_j,p_k$.  The 6 lines not involving a particular index $i$ form
a tetrahedron. The ideal of their union is generated by
the set $F_i$ of 4 cubic forms made by taking
products, 3 at a time, of the $l_{stu}$ with $i \not\in \{s,t,u\}$,
so we may take $d=3$ in the argument above.
\Box

Theorem 1.3 below will show that good partitions
do exist for monomial ideals, and this is the basis of our
method.  The reason that they exist is essentially that a
monomial ideal of codimension $c$ is always contained in an
ideal generated by $c$ elements --- in fact, by
$c$ of the variables.

\vskip .3cm

\centerline{\bf   A method for finding a sparse system of parameters}

\vskip .1cm

\noindent The following theorem is our first main result.
We fix any term order ``$\prec$'' on the polynomial ring $S$
and we write $\,in_\prec({\cal F})\,$ for the set of
initial terms of the polynomials in ${\cal F}$.

\proclaim Theorem 1.3.
Suppose that $in_\prec({\cal F})$ generates an ideal of
codimension $c$. There exist partitions
$\,{\cal F} = {\cal F}_1 \cup \ldots \cup {\cal F}_c \,$
such that for each $i$ the monomials
$in_\prec ({\cal F}_i)$ have a variable in common.
If $\cal F$ is any such partition,
then for almost all $r_{i,f}\in k,$ the
polynomials $(*)$
generate an ideal of codimension $c$.  Further, each
$\,f \in {\cal F}\,$ may be multiplied by any factor
of $in_\prec(f)$ without spoiling this property.

It is known that the hypothesis of Theorem 1.3 is satisfied
if ${\cal F}$ is a Gr\"obner basis of an ideal of codimension $c$.
(See e.g.~[Kalkbrener-Sturmfels 1993]). This suggests the
following algorithm for finding a sparse system of parameters.

 \vskip .3cm

\noindent {\bf Algorithm 1.3'.} \hfill \break
{\sl Input : } A set of generators ${\cal F}$ for an ideal
$I$ of codimension $c$
\hfill \break
Enlarge ${\cal F}$ step by step toward a Gr\"obner basis,  using the
Buchberger algorithm, until $in_\prec({\cal F})$ generates
an ideal of codimension $c$.  Next replace this partial Gr\"obner basis
by a minimal subset which has the same property.
Partition this new ${\cal F}$ into subsets ${\cal F}_i$ as
in  Theorem 1.3,  for example as follows:
Choose  a  prime $(x_{i_1},\ldots, x_{i_c})$ containing $in_{\prec} (\cal F)$.
Such primes exist because every associated prime of a monomial ideal
is generated by a subset of the variables. For $p = 1,2,\ldots,c \,$
define ${\cal F}_p$ inductively to be the set of all elements of
${\cal F} - \cup_{j\leq p}{\cal F}_j$ whose initial terms are
divisible by $x_{i_p}$.

If the polynomials in ${\cal F}$ are homogeneous,
and a homogeneous system of parameters is desired, let $d_i$ be the
maximal degree in ${\cal F}_i$, and multiply each
polynomial in ${\cal F}_i$ by a power of one of the variables in its
own initial term to bring it up to degree $d_i$.

Choose random elements $r_{i,f}$ in $k$, and verify that the
polynomials $f_i$ in $(*)$ generate an ideal of
codimension $c$.  If they do not, try a new random choice.
\hfill \break
\noindent {\sl Output : } The sequence $f_1,\ldots,f_c $. \ \ \Box

\vskip .2cm

Before starting Algorithm $1.3'$, it may be worthwhile to
change to the order for which the initial forms of
${\cal F}$ generate an ideal of largest possible codimension.
This can be done
using the polyhedral methods in [Gritzmann-Sturmfels 1993].

One subtask to be solved in Algorithm $1.3'$ is to
find the minimal prime $(x_{i_1},\ldots, x_{i_c})$.
This is often an easy task,
but we remark that in general it amounts
to solving a combinatorial problem which is NP-complete.
To see this, consider the case where
$in_\prec ({\cal F})$ consists of square-free quadratic
monomials $x_i x_j $, representing the edges of a graph
$G$ with vertex set $\{x_1,\ldots,x_m\}$. A subset $S$
of the vertices of $G$ is called {\it stable} if no two
elements in $S$ are connected by an edge in $G$.
Our subtask amounts to finding a maximal stable set of $G$,
a problem which is known to be NP-complete.

\vskip .3cm

Our proof of Theorem 1.3 is based on the following criterion for a
sequence of sets of polynomials to be what one might describe as a
``setwise system of parameters'':

\proclaim Proposition 1.4.
Let $\, {\cal F}_1,\dots,{\cal F}_c \subset S \,$ be sets of polynomials.
If for every $U \subseteq \{1,\dots,c\}$ the set of polynomials $\,
\bigcup_{j\in U} {\cal F}_j \,$ generates an ideal of codimension
$\geq$ $card\,\, U$, then
for almost every choice of $ \,r_{i,f}\,$ in $\,k\,$  the
polynomials $f_1,\dots,f_c$ in $(*)$ generate an ideal of codimension $c$
in $S$.

\noindent {\bf Remarks. }
\item{(i)} The term ``for almost every choice'' in
Proposition 1.4  (and the term ``random'' in Algorithm 1.1) means
that $r_{i,f}$ can be chosen in some Zariski-open subset of
the coefficient space.
\item{(ii)} Example 1.1 shows that we cannot weaken the hypothesis
by restricting $U$ to be an initial subset of $\{1,\dots,c\}$.
\item{(iii)} The converse of Proposition 1.4 holds after
localizing at any prime containing all the ${\cal F}_i$'s.
(Reason: In a local domain with saturated chain condition,
$codim\, (f_1,\ldots,f_c) = c$ implies
$codim\, (f_1,\ldots,f_d) = d$ for all $d < c $.)
It follows that it holds, even without localizing,
if each ${\cal F}_i$ consists of homogeneous polynomials of
positive degree.
(Reason: The local codimension of the ideal generated by
$ \,\bigcup_{j\in U} {\cal F}_j \,$ is minimized in the
localization at the origin.)

\vskip .2cm

The following example shows that the converse
is false in general. We are grateful to Alicia Dickenstein
for pointing it out.

\noindent {\bf Cautionary Example 1.5. }
Let $m = 3$, ${\cal F}_1 = \{ f_1 = (x_1 + 1) x_2 \} $,
${\cal F}_2 = \{ f_2 = (x_1 + 1) x_3 \} $, and
${\cal F}_3 = \{ f_3 =  x_1 \} $.
Then $f_1,f_2,f_3$ generate an ideal of codimension $3$
but $f_1,f_2$ generate an ideal of codimension $1$.

\vskip .3cm
\centerline{\bf   Proofs}
\vskip .1cm

\noindent {\sl Proof of Proposition 1.4.}
We first show that the conclusion holds in the ``generic situation''.
Let $k'$ be the polynomial ring with
variables $r_{i,f}$ for $f\in {\cal F}_i$ and $i=1,\dots,c$.
We will show that the polynomials $f_1,\dots,f_c$ in $(*)$
generate an ideal of codimension $c$ in $ S \otimes_k k'$.
Equivalently, let $A$ be the affine space with the $r_{i,f}$
as coordinates, and let $B$ be the affine space
with coordinates $x_1,\dots,x_m$. Let $X$ be the
subvariety of $A\times B$ defined by $f_1,\dots,f_c$.  We
will show that the codimension of $X$ is $c$.

Let  $\pi_2 : X \rightarrow B$ be the projection to the second factor
of the product.  The fiber of $\pi_2$ over a point $p\in B$ is a linear space.
Since the $f_i$ involve disjoint sets of
variables $r_{i,f}$, its codimension equals $c$ minus the number of indices
$i $ such that the polynomials in ${\cal F}_i$
all vanish at $p$.  For each subset $U \subseteq \{1,\ldots,c\}$ let
$Y_U$ denote the set of $p\in B$ such that the polynomials in
$ \cup_{j\in U} {\cal F}_j $ all vanish at $p$  but
some polynomial in each ${\cal F}_j$ for $j\notin U$ does not
vanish at $p$. The constructible sets $Y_U$ define a stratification
of $B$ such that over each stratum the fibers of $\pi_2$ have constant
codimension $ c - card\, U$. Therefore the codimension of $X$ in
$A\times B \,$ equals  $\,\,min \bigl\{ \,c - card\, U +codim\, Y_U \,\,: \,\,
U  \subseteq \{1,\ldots,c\} \,\bigr\}$. The
hypothesis in (c) states that if $Y_U$ is nonempty, then the
codimension of $Y_U$ is $\geq card\,\,U$. We conclude that the
codimension of $X$ is $\geq c$.

To conclude the proof, consider the projection  $\pi_1 : X
\rightarrow A$ to the first factor of the product.
We must show that the codimension in $B$ of almost every
fiber of $\pi_1$ is $\geq c$. For any dominant map of
irreducible varieties, almost every fiber has dimension equal
to the dimension of the source minus the dimension of the target.
Thus in our situation almost every fiber has codimension in $B$
equal to the codimension in $A \times B$ of the union of
those components of $X$ that dominate $A$. This second codimension is $\geq$
the codimension of $X$ in $A\times B$, which we have shown to be $\geq c$. \Box

\vskip .2cm

\noindent {\sl Proof of Theorem 1.3. }
The codimension of the ideal generated by
the initial terms of a set of polynomials is $\leq$ that of the ideal
of the polynomials themselves. Thus if the sets  $\,in_\prec({\cal F}_i)\,$
satisfy the hypothesis of Proposition 1.4, then so do the sets ${\cal F}_i$.
Consequently it suffices to treat the case where ${\cal F}$
consists of monomials.  We must find a partition satisfying the
first condition of Theorem 1.3, and show that any such partition
satisfies the hypothesis of Proposition 1.4.

Since every associated prime of a monomial
ideal is generated by a subset of the variables, we may assume
(after renumbering variables if necessary) that
${\cal F}$ is contained  in the ideal $(x_1,\dots,x_c)$.  Since a
monomial is contained in this ideal if and only if
it is divisible by one of the
variables $(x_1,\dots,x_c)$, we  may partition ${\cal F}$ into
subsets ${\cal F}_i$  consisting only of monomials divisible by $x_i$.
The following lemma concludes the proof:

\proclaim Lemma 1.6.  Let $\,{\cal F}_i \,\subseteq \,(x_i) \,\subset S $,
$\, i=1,\dots,c \,$ be sets of monomials.
If $\,{\cal F} \,= \,\cup_{j=1}^c {\cal F}_j\,$ generates an
ideal of codimension $c$, then the ${\cal F}_j$ satisfy the hypothesis
of Proposition 1.4.

\noindent  {\sl Proof. }
For any subset $U$ as in condition (c), let $I_U$ denote the
ideal generated by $\,\cup_{j\in U}\,{\cal F}_j$. The
ideal $I$ generated by ${\cal F}$ is contained in $\, I_U \,
+ \,(\{x_i\}_{i\notin U})$. Since $I$ has codimension $c$, the Principal
Ideal Theorem implies that $I_U$ has codimension $\geq card\,\,U$
as required. \Box

\vskip .3cm

\noindent {\bf Cautionary Example 1.7. }
Let ${\cal F}_i$ and  $\,in_\prec({\cal F}_i)\,$ be as in Theorem 1.3,
and choose coefficients $r_{i,f}$ such that the
linear combinations of initial terms
$$ \sum_{f\in {\cal F}_1} \,r_{1,f}\cdot in_\prec(f) \, ,\quad \ldots
\quad \ldots \, ,\quad \sum_{f\in {\cal F}_c} \, r_{c,f}\cdot in_\prec(f) $$
form a system of parameters. It is tempting to hope that the
$f_i$ in $(*)$, made with  the \underbar{same} coefficients $r_{i,f}$,
would  also form a system of parameters.
This is  not true:  If ${\cal F}_1 = \{x_1^2- x_2^2\}$
and  ${\cal F}_2 = \{x_1 x_2, x_2^2\}$  then $x_1^2, x_1 x_2+ x_2^2$
is a regular sequence but $x_1^2- x_2^2,x_1 x_2+ x_2^2$ is not. \Box

\vskip .2cm

We close Section 1 with two propositions showing that our examples
1.1 and 1.2 are minimal in a certain sense.
The proofs are straightforward and will be omitted.

\proclaim Proposition 1.8.  Suppose that ${\cal F} \subset S$
generates an ideal of codimension $c$ and that
${\cal F}$ cannot be partitioned into subsets ${\cal F}_1,\dots,{\cal F}_c$
such that for some choice of coefficients
$r_{i,f}$ the polynomials $(*)$ form a system of parameters.
\item{(a)} The set ${\cal F}$  can be replaced by a set of $\,c+1\,$ linear
combinations of the  elements of ${\cal F}$ having the same property,
possibly after reducing $c$.
\item {(b)} Factoring out $m-c$ general
linear forms, the number of variables of $S$ may be taken to be $\, c $.

Thus the critical case concerns sets of $c+1$
polynomials generating an ideal of codimension $c$ in $c$
variables.  It is most interesting to look at the case of
homogeneous polynomials.  Example 1.2 is of this kind,
with $c=4$, but there are no such examples with $c\leq 3$:

\proclaim Proposition 1.9.  If ${\cal F} = \{f_1,f_2,f_3,f_4\}$ is a set of
 homogeneous polynomials in 3 variables, generating an ideal of
codimension 3 in $S=k[x_1,x_2,,x_3]$,
then there is a partition of ${\cal F}$ into 3 subsets ${\cal F}_i$
such that the polynomials in $(*)$ form a system of parameters.

\beginsection 2.   Systems of parameters modulo an ideal

We now turn to the general case of our problem, keeping the
notation $S := k[x_1,\ldots,x_m]$ as in the introduction.  Let
$J\subset S$ be an ideal, let $R = S/J$ and let ${\cal F} \subset S$
be a finite subset.  Let $I$ be the ideal generated by $\cal F$,
and suppose that the codimension of $I$ modulo $J$ is $c$ in the
sense that $c=codim (I+J) - codim (J)$.

We say that $f_1,\ldots, f_c \in I$ is a maximal system
of parameters for $I$ modulo $J$ if
$\,codim(P+(f_1,\ldots, f_c)) \,\geq\,  codim (I+J)\,$ for
every minimal prime $P$ of $J$.   (This notion is most natural
if the ideal J is unmixed.) We wish to choose as sparse as
possible a maximal system of parameters for $I$ modulo $J$.

The simplest and most common problem calling for systems of
parameters is that of finding a Noether normalization
for a homogeneous ideal. If $c$
is the Krull dimension of $R$, this is the problem
of finding elements $f_1,\ldots ,f_c$ in $R$ such that
$R$ is a finitely generated module over the subring
$k[f_1,\dots, f_c] \subseteq R $. (The elements $f_1,\ldots ,f_c$
are then necessarily algebraically independent, so that
the subring is isomorphic to a polynomial ring.
See [Logar 1988] and [Dickenstein, Fitchas, Giusti, and Sessa 1991]
for a discussion
from a computer algebra point of view.) We will focus primarily on this
case, but first we present a method for handling the general
problem. The approach
differs from that of Section 1 in
that it chooses one $f_i$ at a time, essentially employing
overlapping sets ${\cal F}_i$.

\proclaim Greedy Algorithm 2.1. \hfill \break
{\sl \ \quad Input : } \rm
A set of generators ${\cal F}$ for an ideal $I$ of
codimension $c$ modulo $J$.
\hfill \break
\rm
Let ${\cal F}_1$ be a minimal subset of $\cal F$ such that
${\cal F}_1$ is not contained in any minimal prime of $J$
of maximal dimension.
Let $f_1$ be a sufficiently general
combination of the elements of ${\cal F}_1$ so that
the codimension of $J+(f_1)$ is larger than that of $J$.  Let ${\cal F}'$ be
the result of dropping any one of the elements of ${\cal F}_1$
from ${\cal F}$.  Replace $J$ by $J+(f_1)$,
 replace ${\cal F}$ by
${\cal F}'$, and iterate the process. \hfill \break
{\sl Output : } The sequence $f_1,\ldots,f_c$.
\vskip .1cm

In case the data ${\cal F}, I, J$ are homogeneous
and the output desired is homogeneous,
but not all the polynomials of ${\cal F}$ are of the same degree,
the ``sufficiently general linear combination''
would have to have coefficients that are polynomials of
varying degrees.  The following variation may be an improvement:

Given a minimal subset ${\cal F}_1\subseteq {\cal F}$
not contained in any minimal prime of $J$
of maximal dimension, replace it with a set of elements whose initial
forms are not contained in any minimal prime of $in(J)$ of maximal
dimension.  (This may be done by moving step by step toward a
Gr\"obner basis of $J+({\cal F}_1)$, using the Buchberger algorithm,
until the codimension of the initial ideal is larger than that of the
initial ideal of $J$.)  Then, if homogeneous output is desired, each
element of ${\cal F}_1$ not of maximal degree may be multiplied
by variables dividing its initial term to bring all the elements of
${\cal F}_1$ to the same degree before forming the linear combination
as above.

We now turn to the special case of Noether normalization.
The following well-known
version of Hilbert's Nullstellensatz makes clear the nature
of our task:

\proclaim Proposition 2.2.
Let $J$ be a homogeneous ideal of $S$, and set $R = S/J$. Let
$X \subset {\rm \bf P}^{m-1}$
be the corresponding projective algebraic set.  Suppose that
the ground field $k$ is algebraically closed.
If $f_1,\ldots ,f_c \in R$ are homogeneous polynomials, then
$R$ is a finitely generated module over the  subring
$k[f_1,\dots, f_c] \subseteq R $ if and only if the system of equations
$\,f_1({\bf x}) = \ldots = f_c ({\bf x}) = 0 \,$
has no solution in $X$.

\noindent {\sl Proof Sketch: }
By the Nullstellensatz the condition that there are no solutions
is equivalent to the condition that $R/(f_1,\dots, f_c )$ is a
finite dimensional vector space.  Because R is graded and
zero in negative degree,
a basis for this
space lifts to a finite set of generators for $R$ over the subring
$k[f_1,\dots, f_c]$.  \Box

In the Noether normalization problem one usually wants
the $f_i$ to be linear forms.  We will henceforth
consider only this case, and suppose that
${\cal F} = \{ x_1,\ldots,x_m \}$, so that $I = ({\cal F})$
is the irrelevant ideal. In this situation Algorithm 2.1
has the effect of reducing at each step the number of
variables to be considered, and this increases its efficiency.
Here is a monomial example, which also suggests a possibility for
improving the algorithm:

\vskip .2cm

\noindent {\bf Example 2.3. }
Let $m=6$, $c=2$, ${\cal F} = \{x_1,x_2,x_3,x_4,x_5,x_6\}$, and
$$ \eqalign{
\,J \,\,\, = \quad & \bigl(\, x_1 x_2, \,x_1 x_3,\, x_2 x_3,
\,x_2 x_4, \,x_2 x_5, \,x_3 x_4, \,x_3 x_5,
\,x_4 x_5 , \,x_4 x_6 , \,x_5 x_6 ) \quad = \cr
 (x_1,x_2,&x_4,x_5) \,\cap\,
(x_1,x_3,x_4,x_5) \,\cap\,
(x_2,x_3,x_4,x_5) \,\cap\,
(x_2,x_3,x_4,x_6) \,\cap\,
(x_2,x_3,x_5,x_6) .\cr}$$
In the first step the unique optimal choice is ${\cal F}_1 = \{x_1, x_6\}$.
We set $\,f_1 \, = \,x_1 - x_6 $ and remove
$x_1$ from ${\cal F}$. Repeating the procedure with
$$ J + (f_1) \quad = \quad \bigl( \,x_1 - x_6, x_2 x_3,  x_2 x_4,  x_2 x_5,
 x_2 x_6,  x_3 x_4,  x_3 x_5, x_3 x_6,  x_4 x_5, x_4 x_6, x_5 x_6 \,\bigr),$$
we must use all remaining variables: $\,{\cal F}_2 = \{ x_2,x_3,x_4,x_5,x_6\}$.
For the second parameter we can choose,
for instance, $\,f_2 \,=\,x_2 + x_3 + x_4 + x_5 + x_6 $.\ \ \Box

\vskip .2cm

In general Algorithm 2.1  is ``too greedy'': it does not perform
optimally with respect to sparseness. For example,
$$ \hat f_1 \,= \,x_1 + x_2 + x_3,\quad  \hat f_2 \,= \,x_4 + x_5 + x_6 $$
is also a Noether normalization for the ideal $J$ in Example 2.3.
It has a total of only six non-zero terms, while the output
$f_1,f_2$ of Algorithm 2.1 has seven non-zero terms.

This shows that the subtlety of sparse Noether normalization
is not completely captured by Algorithm 2.1.
The remainder of this paper is devoted to a
more thorough combinatorial analysis, leading to an optimal result.

\vskip .4cm

\centerline{\bf Sparsity of Linear Subspaces}

\vskip .1cm
\noindent
We begin with some remarks on the
notion of sparseness itself.
{}From Proposition 2.2 we see that the problem of
choosing a space of linear forms $f_1,\ldots, f_c \in
I$ of given degree $d$ that are a homogeneous system
of parameters modulo a homogeneous
ideal $J$, is equivalent to the following geometric
problem: given an algebraic set $X$  of dimension
$c-1$ in a projective space ${\rm\bf P}^{m-1}$, find a
linear subspace $L$ of codimension $c$ not meeting
$X$.  Equivalently, we may think of $L$ as
coming from a linear subspace $M$ of
an affine space ${\rm\bf A}^{m}$, which is
supposed to meet the cone over $X$ only in the origin.

We wish to choose $M$ to be as sparse as possible,
relative to some given system of coordinates for
${\rm\bf P}^{m-1}$.  There are several
possible  definitions of sparseness, and they
conflict with one another.  In general, if we agree
on a way to  represent the space $M$, then we can
speak of a space allowing the sparsest possible
representation in this  form. Perhaps the three most
obvious representations are these: $M$ might be
represented by the coordinates of a basis of
 $M$ ({\it basis representation}), by the
coordinates of a basis for the space
$M^{\perp}$ of linear
functionals  vanishing on $M$
({\it cobasis representation}),
or by Pl\"ucker
coordinates, the maximal  minors of some matrix
representing the basis or dual basis
({\it Pl\"ucker representation}). In each case the
number of nonzero coordinates is a measure of
sparseness --- we call them {\it basis sparseness},
{\it cobasis sparseness}, and {\it Pl\"ucker sparseness}
respectively.

It is not hard to show that
for 1-dimensional subspaces (and thus also
for hyperplanes) the three measures
of sparseness agree in the sense that all three
choose the same space as the sparsest in a particular
set of subspaces.
But in general no two of these measures agree on
naming the sparsest subspace, as may be seen
from the following examples.  In each case the
subspace considered is the row space of the
given matrix. As we will not make use of these
facts, we leave their verification to
the interested reader.

First, to show that Pl\"ucker
sparseness does not agree with basis
sparseness:
The space $M_1$ represented by the matrix
$$
M_1 \leftrightarrow
\pmatrix{ 1 & 0 & 0 & 0 & 0 & 0 & 0\cr
             0 & 1 & 1 & 1 & 1 & 1 & 1 \cr}
$$
has basis sparseness 7 and Pl\"ucker sparseness
6, while the space $M_2$ represented by
$$
M_2 \leftrightarrow
\pmatrix{ 1 & 1 & 1 & 0 & 0 & 0 & 0\cr
             0 & 0 & 0 & 1 & 1 & 1 & 0 \cr}
$$
has basis sparseness 6 and Pl\"ucker sparseness
9.

        The sparseness of a space $M$ is the same as
the cobasis sparseness of $M^{\perp}$, so the spaces
$M_1^{\perp}$ and $M_2^{\perp}$ illustrate the
same point for cobasis sparseness.

It is harder to give examples in which basis
sparseness and cobasis sparseness disagree, but
the reader may check that if $L_1$ and $L_2$ are
the 3-dimensional subspaces of a 9-dimensional
vectorspace $V$
represented by the matrices

$$
L_1 \leftrightarrow
\pmatrix{ 1 & 1 & 1 & 1 & 1 & 1 & 1 & 1 & 1\cr
          0 & 0 & 0 & 1 & 1 & 1 & 2 & 2 & 2\cr
                       0 & 1 & 2 & 0 & 1 & 2 & 10 & 11 & 12\cr}
$$
and
$$
L_2 \leftrightarrow
\pmatrix{ 1 & 1 & 1 & 1 & 1 & 1 & 1 & 1 & 1\cr
          0 & 0 & 0 & 1 & 1 & 1 & 2 & 3 & 4\cr
                       0 & 1 & 2 & 0 & 1 & 2 & 10 & 11 & 12\cr}
$$
then the basis sparseness of $L_1$ is $6+6+7 = 19$,
whereas that of $L_2$ is $6+6+6=18$.  The cobasis
sparseness of each is $3+3+3+4+4+4 = 21$.
Now consider the subspaces
$$
M_3 = L_1\oplus L_2^{\perp}\quad\quad{\rm and}\quad\quad
M_4 = L_2\oplus L_1^{\perp} \quad \,\,\, \hbox{ in } \quad
V\oplus V^*.
$$
The basis sparsenesses of these
spaces are $40 = 19+21$ and $39=18+21$ respectively.
But as
$M_3^{\perp}=(L_1\oplus L_2^{\perp})^{\perp} =
L_1^{\perp}\oplus L_2$, and similarly for $M_4$,
the cobasis sparsenesses for $M_3$ and $M_4$ are
$39=18+21$ and $40 = 19+21$,
reversing the order of the basis sparsenesses.

\vskip .4cm

\centerline{\bf Sparse Noether Normalization using Chow Forms}

\vskip .1cm

\noindent
We now return to the problem of Noether Normalization.
We will work in terms of basis sparseness of the
space generated by the linear forms in the solution to
the Normalization problem; our discussion can be
adapted, by considering different expressions of the Chow form,
to cobasis or Pl\"ucker sparseness as well.
Let $J$ be a homogeneous unmixed ideal in $S$,
and let $X$ denote its projective variety in ${\bf P}^{m-1}$.
Changing notation somewhat, we suppose that
$X$ has degree $p$ and dimension $d-1$.
We will show how to
compute a Noether normalization consisting of linear forms
$$ f_i \,\, := \,\,c_{i1} x_1 \, + \,c_{i2} x_2 \, + \,\ldots \, +
c_{im} x_m \, , \qquad i = 1,2,\ldots,d \eqno (2.1) $$
 which is
optimal in the sense that the number of
non-zero coefficients $c_{ij}$ is minimal --- that is,
the basis sparseness of the space spanned by the $f_i$ is
minimal. We call the number of nonzero
$c_{ij}$ the {\it Noether complexity} of $X$.

Let $\,R_X = R_X( c_{ij}) = R_X (f_1,\ldots,f_d)$
denote the {\it Chow form} of $X$. This
classical polynomial is characterized by the property that it
vanishes if and only if the linear subspace defined by
$\,f_1({\bf x}) = \ldots = f_d ({\bf x}) = 0 \,$
meets $X$; see e.g.~[Shafarevich 1977, \S I.I.6.1]
and the references given in [Caniglia 1990].
The notation concerning Chow forms tends to vary from
author to author. The specific notation to be employed here
is taken from [Kapranov-Sturmfels-Zelevinsky 1992] and
[Sturmfels 1993], namely, we
express $R_X$ as a polynomial in {\it brackets}
$\, [\,i_1 \,i_2 \,\ldots \,i_d\,]$, $\,1 \leq i_1 <
\ldots < i_d \leq m $. These are the $d \times d$-minors
of the $d \times m $-matrix $(c_{ij})$, or, equivalently,
the Pl\"ucker coordinates of the codimension $k$ flat defined by (2.1).
By Proposition 2.2, the Noether normalization problem for $R$ is equivalent
to the problem of finding a non-root of the Chow form $R_X$.

\vskip .2cm

\noindent {\bf Example 2.4. }
{\sl (Hypersurfaces, $d = m-1 $)} \hfill \break
Suppose $J$ is the principal ideal generated by a homogeneous
polynomial $\, F(x_1,x_2,\ldots,x_m )$, defining a hypersurface
$X \subset {\bf P}^{m-1}$. In terms of brackets its Chow form equals
$$ R_X \quad = \quad
F  \bigl( \, [2 3 4 \ldots m] \, , \,
- [1 3 4 \ldots m] \, , \,[1 2 4 \ldots m] \, , \,\,\ldots \,, \,\,
(-1)^{m-1} [2 3 \ldots m - 1] \, \bigl) .  \eqno (2.2) $$
The Noether normalizations of $X$ are precisely the
$(m-1) \times m $-matrices $(c_{ij})$ for which the
bracket polynomial (2.2) does not vanish.
\ \Box

\vskip .2cm

It is our objective to compute a $ d \times m $-matrix
$(c_{ij})$, which is a non-root of the Chow form $R_X$,
and which is as sparse as possible with this property.
Since $(c_{ij})$ must have maximal rank $d$,
the Noether complexity of $X$ is at least $d = dim(X)+1$.
It is exactly $d$ if and only if $X$ is
in Noether normal position with respect to some coordinate flat.

\proclaim Observation 2.5.
The coordinate forms $\, x_{i_1}, x_{i_2},\ldots, x_{i_d}\,$
are a Noether normalization for $X$ if and only if
the bracket power $ [\,i_1 \,i_2 \,\ldots \,i_d\,]^p$ appears
with non-zero coefficient in $R_X$.

Let $ V = \{ c_{ij} : 1 \leq i \leq d , 1 \leq j \leq m \}$ denote
the set of variables. For any polynomial $f \in k[V]$ we define
a simplicial complex $\Delta(f)$ as follows:
A subset $W \subset V$ is a face of $\Delta(f)$ if and only if
there exists a non-root of $f$ whose zero coordinates are
precisely $W$. Equivalently, $W$ is not a face of $\Delta(f)$
whenever $f$ lies in the ideal generated by $W$. If we write
$supp(m)$ for the set of variables dividing a monomial $m$,
we see that the maximal faces of $\Delta(f)$ are the complements
of the minimal sets of the form $supp(m)$ where $m$ is a monomial
of $f$.  In particular, for each monomial $m$, the complex
$\Delta(m)$ is a simplex, consisting of all subsets of
$\,V \setminus supp(m)$. Thus we get the first statement of the following:

\proclaim Lemma 2.6.
Let $f $ be a homogeneous polynomial in $k[V]$.
Then $\Delta(f)$ is the union of the simplices
$\Delta(m)$, where $m$ ranges
over all monomials of $f$ with minimal support.
Also, $\Delta(f)$ is the union of the simplices
$\Delta(in_\prec (f))$, where $\prec$ ranges
over all term orders on $k[V]$.

\noindent {\sl Proof: }
We have already proved the first statement.  To prove the second
it suffices to observe that because $f$ is homogeneous,
every monomial of $f$ with minimal support is the
initial monomial of $f$ with respect to some term order.
\Box

 \vskip .1cm

This lemma together with the above observations implies
the following theorem.

\proclaim Theorem 2.7. The Noether complexity of a projective
variety $X$ equals the least number of variables $c_{ij}$
appearing in any initial monomial $\,in_\prec (R_X) =
\prod c_{ij}^{\nu_{ij}}$ of its Chow form. \hfill \Box

\vskip .2cm

\noindent {\bf Example 2.3.~(continued)}
The reducible curve $X \subset {\bf P}^5$ defined by
$J$ has Chow form
$$ \eqalign{
 & R_X  \quad = \quad
 [\,1\,4\,]\cdot [\,1\,5\,]\cdot [\,1\,6\,]
\cdot [\,2\,6\,] \cdot [\,3\,6\,] \quad = \cr
&  (c_{11} c_{24} - c_{14} c_{21}) \cdot
(c_{11} c_{25} - c_{15} c_{21}) \cdot
(c_{11} c_{26} - c_{16} c_{21})\cdot
(c_{12} c_{26} - c_{16} c_{22}) \cdot
(c_{13} c_{26} - c_{16} c_{23}).\cr} $$
The coefficient matrices of $f_1,f_2$ and $\hat f_1, \hat f_2$
considered above are seen to be non-roots of $R_X$.
The Noether complexity of the curve $X$ is
equal to six   \ (cf.~Theorem 2.7). \ \Box

The most systematic approach to solving our problem
would be to explicitly compute the Chow form $R_X$.
By the results of [Caniglia 1990], this can be done
in single exponential time (in $m$). Theorem 2.7 implies
that the Noether complexity and an optimal Noether normalization
for $X$ can be computed in single exponential time.

Unfortunately, this approach is not useful in practice,
since the complete expansion of the Chow form into monomials
$\prod c_{ij}^{\nu_{ij}}$ is usually too big.
Hence the Caniglia algorithm is only of theoretical interest
with regard to our problem. In fact,
the problem of computing the Noether complexity
of a monomial ideal is NP-hard. The following
proof of this fact has been pointed out to us by
Jesus DeLoera. Let $G$ be any graph on $V =
\{x_1,\ldots,x_n\}$ and $I_G$ the ideal
generated by all square-free cubic monomials
$x_i x_j x_k$, and all $x_i x_j$ not
corresponding to an edge of $G$ (this is the
Stanley-Reisner ideal of $G$ viewed as a
simplicial complex). The Noether complexity of
$I_G$ equals the minimal number $ 2 n - |S_1| -
|S_2|$, where $S_1,S_2 \subset V$ ranges over
all disjoint pairs of stable sets of $G$. Here
$S_i$ indexes the zero entries in row $i$ of a
sparsest Noether normalization $(c_{ij})$. If we
had a polynomial time algorithm for finding
$S_1$ and $S_2$, then we could solve the NP-hard
problem of computing a maximal stable set in any
graph $G_1$ as follows. Let $G_2$ be a disjoint
copy of $G_1$, and let $G$ be the graph obtained
from their union union $G_1 \, \cup \,G_2$ by
connecting each vertex of $G_1$ with each vertex
of $G_2$. Applying our algorithm to $G$ we obtain
a maximal stable set $S_i$ for $G = G_i$, $i=1,2$.

For practical computations we propose an approach using
(truncated) Gr\"obner basis computations for the ideal $J$.
The following proposition, which is easily derived from the
proof of Lemma 2.6, shows that the Noether complexity of
a homogeneous ideal is bounded above by the Noether complexity
of any of its initial ideals.

\proclaim Proposition 2.8. Let $(c_{ij})$
be any Noether normalization of the initial monomial ideal $in_\omega(J)$,
where $\omega = (\omega_1,\ldots,\omega_n) \in {\bf Z}^{m}$
represents any term order for $J$. Then, for almost all $t \in k$,
the matrix $\,(c_{ij} \cdot t^{\omega_j}) \,$
is a Noether normalization of $J$.

{}From this we get the following lifting algorithm:
Choose any term order on $S$, and
compute a truncated Gr\"obner basis $\{g_1,\ldots,g_r\} $ for $J$,
subject to the truncation condition that the monomial ideal
$\, L = \bigl(in(f_1), \ldots, in(f_r)\bigl) $ has the same radical as
$in_\omega(J)$. Compute the Chow form $R_L$ of $L$, e.g., using
method in Proposition 3.4 of [Sturmfels 1992]. The bracket monomial
$R_L$ has precisely the same bracket factors as $R_{in(J)}$. We have
$$ \Delta( R_L) \quad = \quad
\Delta(R_{in(J)}) \quad \subseteq \quad \Delta(R_X) .\eqno (2.3) $$
Choose any maximal face of the simplicial complex $\Delta(R_M)$.
This gives rise to a Noether normalization for $in(J)$ and,
using Proposition 2.8, we get a Noether normalization for $J$.

In order to find a sparser Noether normalization we may repeat
this procedure for as many different term orders as we can.
In fact, whenever this is feasible, one might like to compute a
{\it universal} Gr\"obner basis  ${\cal U}$ for $J$,
that is, a finite subset of $J$ which is a Gr\"obner basis
simultaneously for all term orders on $S$.
{}From a universal Gr\"obner basis and the knowledge
of the maximal faces of $\Delta(R_M)$ we can read off the minimum
of the Noether complexities of all initial ideals of $J$.
However, this minimum will generally not agree with
the Noether complexity of $J$, as the following example shows.

\vskip .2cm

\noindent {\bf Cautionary Example 2.9.}
A homogeneous ideal $J$ whose Noether complexity
is smaller than the Noether complexity of any of its initial
ideals $in(J)$. Let $m=6$ and consider
$$ J \quad = \quad
( x_2 x_5 - x_1 x_6 , x_3, x_4) \,\,\cap \,\,
( x_1 x_4 - x_3 x_5 , x_2, x_6) \,\,\cap \,\,
( x_3 x_6 - x_2 x_4 , x_1, x_5) . $$
The variety $X$ of $J$ is a union of three toric surfaces in ${\bf P}^5$.
Here the Chow form equals
$$ R_X \quad = \quad
\bigl([126][156] - [125][256]\bigr) \cdot
\bigl([135][345] - [145][134]\bigr)\cdot
\bigl([234][246] - [236][346]\bigr) . $$
The matrix
$$ (c_{ij}) \quad = \quad
\pmatrix{ 1 & 0 & 0 & 1 & 0 & 0 \cr
             0 & 1 & 0 & 0 & 1 & 0 \cr
             0 & 0 & 1 & 0 & 0 & 1 \cr} $$
is a non-root of $R_X$ and hence defines a
Noether normalization. It is optimal because each term
appearing in the complete expansion of $R_X$
contains at least six of the variables $c_{ij}$.
This proves that the Noether complexity of $X$ equals six.

The ideal $J$ has six  distinct initial
ideals, each of which is isomorphic to
$$ in(J) \quad = \quad
( x_2 x_5, x_3, x_4) \,\,\cap \,\,
( x_1 x_4 , x_2, x_6) \,\,\cap \,\,
( x_2 x_4 , x_1, x_5) . $$
The complete expansion of the initial Chow form
$$ R_{in(J)} \quad = \quad
[126][156][135][345][236][346] $$
has $13,452$ terms. Each of these terms contains at least
eight variables. Hence the Noether complexity of
each initial ideal of $J$ equals eight. \Box

\vfill \eject
\vskip .7cm

{\baselineskip=12pt

\centerline {\bf References}

\vskip .2cm \par \frenchspacing

\noindent D.~Bayer and M.~Stillman: Computation of Hilbert functions,
{\sl J.~Symb.~Comput.} {\bf 14} (1992) 31--50.
 \vskip .2cm

\noindent  A.~Bigatti, M.~Caboara, and L.~Robbiano:
On the computation of Hilbert-Poincar\'e series,
{\sl Applicable Algebra in Engineering, Communications
and Computing}, {\bf 2} (1993) 21--33.
\vskip .2cm

\noindent L.~Caniglia: How to compute the Chow form of
an unmixed polynomial ideal in subexponential time,
{\sl Applicable Algebra Eng.~Commun.~Comput.} {\bf 1} (1990) 25--41.
 \vskip .2cm

\noindent D.~Eisenbud, C.~Huneke, and W.~Vasconcelos: Direct
methods for  primary decomposition. {\sl Inventiones
Math.} {\bf 110} (1992) 207--236. \vskip .2cm

\noindent D.~Eisenbud: Open problems in computational algebraic
geometry. To appear in the {\sl Proceedings of the Cortona Conference on
Computational Algebraic Geometry}, ed. D.~Eisenbud and
L.~Robbiano,  Cambridge Univ. Press, Cambridge (1993). \vskip .2cm

\noindent A.~Dickenstein, N.~Fitchas, M.~Giusti, and C.~Sessa:
The membership problem for polynomial ideals is solvable in subexponential
time, {\sl Discrete Applied Math.} {\bf 33} (1991) 73--94.
 \vskip .2cm

\noindent M.~Kalkbrener and B.~Sturmfels:
Initial complexes of prime ideals.
To appear in {\sl Advances in Mathematics} (1993).
 \vskip .2cm

\noindent M.~Kapranov, B.~Sturmfels, and A.~Zelevinsky:
Chow polytopes and general resultants,
{\it Duke Mathematical Journal} {\bf 67} (1992) 189--218. \vskip .2cm

\noindent T.~Krick and A.~Logar: An algorithm for the
computation of the radical of an ideal in the ring of polynomials.
Proceedings 9th AAECC, Lecture Notes in Computer Science  539,
Springer-Verlag, New York (1991) pp.~195-205. \vskip .2cm

\noindent P.~Gritzmann and B.~Sturmfels:
Minkowski addition of polytopes: Computational complexity
and applications to Gr\"obner bases.
{\sl SIAM J.~Discr.~Math.} {\bf 6} (1993) 246-269.
\vskip .2cm

\noindent
A.~Logar: A computational proof of the Noether Normalization lemma,
{\sl Springer Lecture Notes in Computer Science},
{\bf 357}, Proceedings AAECC-6 Rome 1988, pp.~259-273. \vskip .2cm

\noindent  I.~Shafarevich: Basic Algebraic Geometry,
Springer Verlag, New York, 1977.
\vskip .2cm

\noindent B.~Sturmfels: Sparse elimination theory.
To appear in the {\sl Proceedings of the Cortona Conference on
Computational Algebraic Geometry}, ed. D.~Eisenbud and
L.~Robbiano,  Cambridge Univ. Press, Cambridge (1993). \vskip .2cm

\noindent W.~Vasconcelos: Constructions in commutative algebra.
To appear in the {\sl Proceedings of the Cortona Conference on
Computational Algebraic Geometry}, ed. D.~Eisenbud and
L.~Robbiano,  Cambridge Univ. Press, Cambridge (1993). \vskip .2cm

\noindent W.~Vasconcelos: The top of a system of equations. To
appear in {\sl Boletin Soc.~Mex.~Matematica } (issue
dedicated to Jose Adem) (1993)
\vskip .2cm

}
\bye